\documentstyle[amssymb,makeidx,12pt]{article}

\newtheorem{th}{Theorem}[section]

\newtheorem{defn}[th]{Definition}
\newenvironment{defn-new}{\begin{defn} \em}{\end{defn}}
\newtheorem{rem}[th]{Remark}
\newenvironment{rem-new}{\begin{rem} \em}{\end{rem}}
\newtheorem{ex}[th]{Example}
\newenvironment{ex-new}{\begin{ex} \em}{\end{ex}}
\newtheorem{exer}[th]{Exercise}
\newenvironment{exer-new}{\begin{exer} \em}{\end{exer}}
\newtheorem{agr}[th]{Agreement}
\newenvironment{agr-new}{\begin{agr} \em}{\end{agr}}
\newtheorem{pbm}[th]{Problem}
\newenvironment{pbm-new}{\begin{pbm} \em}{\end{pbm}}

\makeatletter \@addtoreset{equation}{section} \makeatother

\begin{document}

\begin{center}
$(1,2)${\bf -NULL BERTRAND CURVES IN MINKOWSKI SPACETIME}

\bigskip

{\bf Mehmet G\"{o}\c{c}men and Sad\i k Kele\c{s}}

\bigskip
\end{center}

\noindent {\bf Abstract}. In this paper we study null Bertrand curves in $%
R_{1}^{4}$ under the assumption the curve has a Cartan frame. We show that
if the derivative vectors of the null Cartan curve in $R_{1}^{4}$ is
linearly independent, then this curve is not a Bertrand curve. Since then
the already known notion of \ null Bertrand curves in $R_{1}^{4}$ occurs
only if the derivative vectors of the curve is linearly dependent. We will
introduce an idea of Bertrand curves and abiding by this idea we bring to
light under which conditions a null Cartan curve in $R_{1}^{4}$ is a
Bertrand curve.



\section{Introduction\label{sect-intro}}

Hiroo Matsuda ve Shinsuke Yorozu in their paper \cite{matsuda} introduced a
new type of curves called special Frenet curves and proved that a special
Frenet curve in $R^{n}$ is not a Bertrand curve if $n\geqslant 4$ . They
also improved an idea of generalized Bertrand curve in $R^{4}$. Particularly
they characterized $\left( 1,3\right) $-Bertrand curves in $R^{4}$ and
illustrate this type of curves with an example. Honda-Inoguchi \cite{honda}
and Inoguchi-Lee \cite{lee} have done some work on a pair of null curves $%
\left( C,\overline{C}\right) $, called a null Bertrand pair and their
relation with null helices in $R_{1}^{3}$ . For a null curve in the
Minkowski n-space, Ferrandez, Gimenez and Lucas in \cite{angel} assumed the
linear independence of the derivative vectors of the curve to get a unique
Cartan frame for the curve. On the other hand in \cite{sasaki} Makoto Sakaki
uttered that the assumption in the Theorem 2 in \cite{angel} can be lessened
for a null curve to have a unique Cartan frame in $R_{1}^{n}$. In \cite%
{ceylan} A.Ceylan \c{C}\"{o}ken and \"{U}nver \c{C}ift\c{c}i studied null
Bertrand Curves in Minkowski Spacetime.

Here we show that a null Cartan curve in $R_{1}^{4}$ is not a Bertrand
curve\ if the derivative vectors of the curve is linearly independent and
give a characterization of $\left( 1,2\right) $-Bertrand curves in $%
R_{1}^{4}.$

\section{Preliminaries\label{pre}}

Let $c:I\longrightarrow R_{1}^{n}$, $n=m+2$ , be a null curve parametrized
by the pseudo-arc parameter such that $\left\{ c^{\prime }\left( s\right)
,c^{\prime \prime }\left( s\right) ,...,c^{\left( n\right) }\right\} $ is a
basis of $T_{c\left( s\right) }R_{1}^{n}$ for all $s$. Then there exist only
one Frenet frame satisfying the equations

\begin{eqnarray}
L^{\prime } &=&W_{1},  \nonumber \\
N^{\prime } &=&k_{1}W_{1}+k_{2}W_{2},  \nonumber \\
W_{1}^{\prime } &=&-k_{1}L-N,  \label{1} \\
W_{2}^{\prime } &=&-k_{2}L+k_{3}W_{3},  \nonumber \\
W_{i}^{\prime } &=&-k_{i}W_{i-1}+k_{i+1}W_{i+1}\ \ \ i\in \left\{
3,...,m-1\right\}  \nonumber \\
W_{m}^{\prime } &=&-k_{m}W_{m-1},  \nonumber
\end{eqnarray}%
and verifying

\noindent (i) For $2\leqslant i\leqslant m-1$, $\left\{ c^{\prime
},c^{\prime \prime },...,c^{\left( i+2\right) }\right\} $ and $\left\{
L,N,W_{1},...,W_{i}\right\} $ have the same orientation.

\noindent (ii) $\left\{ L,N,W_{1},...,W_{m}\right\} $ is positively oriented
\cite{angel}.

Assume that $c:I\rightarrow R_{1}^{n}$ be a null Cartan curve. Then for the
Cartan curvatures of $C$, we obtain $k_{2}<0$, $k_{i}>0$ for all $i$ $\in
\left\{ 3,...,m-1\right\} $, and $k_{m}>0$ or $k_{m}<0$ according to $%
\left\{ c^{\prime },c^{\prime \prime },c^{\left( 3\right) },...,c^{\left(
n\right) }\right\} $ is positively or negatively oriented, respectively \cite%
{angel}.

Let $\left( C,\overline{C}\right) $ be a pair of \ framed null Cartan curves
in $R_{1}^{4}$, with distinguished parameters $s$ and $\overline{s}$
respectively. This pair is said to be a null Bertrand pair if their
principal normal vector fields are linearly dependent. The curve $\overline{C%
}$ is called a Bertrand mate of $C$ and vice versa. A framed null curve is
said to be a null Bertrand curve if it admits a Bertrand mate. By this
definition there exist a functional relation $\overline{s}=\varphi \left(
s\right) $ for a null Betrand pair $\left( C,\overline{C}\right) $ such that
$\overline{W}\left( \overline{s}\right) =\mp W\left( s\right) $, i.e., the
normal lines coincide at their corresponding points \cite{duggal}.

\section{Bertrand curves in $R_{1}^{4}$}

Let $C$ be a null Cartan curve in $R_{1}^{4}$. After a straightforward
computation we can write the following equation

\[
\left(
\begin{array}{c}
c^{\prime } \\
c^{\prime \prime } \\
c^{\left( 3\right) } \\
c^{\left( 4\right) }%
\end{array}%
\right) =\left(
\begin{array}{cccc}
1 & 0 & 0 & 0 \\
0 & 0 & 1 & 0 \\
-k_{1} & -1 & 0 & 0 \\
-k_{1}^{%
{\acute{}}%
} & 0 & -2k_{1} & -k_{2}%
\end{array}%
\right) \left(
\begin{array}{c}
L \\
N \\
W_{1} \\
W_{2}%
\end{array}%
\right)
\]%
\ \ If we denote the above matrix by $P$ , then $\left\vert
P\right\vert =-k_{2}$. Here we choose the spacelike vector $W_{2}$
uniquely so that the pseudo-orthonormal basis $\left\{
L,N,W_{1},W_{2}\right\} $ is positively oriented. If we assume that
$\left\{ c^{\left( i\right) }\right\} _{1\leq i\leq 4}$ is
positively oriented, then we get $k_{2}<0$. Notice that
if $k_{2}=0$, then $\left\vert P\right\vert =0$ and this says that the set $%
\left\{ c^{\left( i\right) }\right\} _{1\leq i\leq 4}$ is linearly dependent.

\begin{th}
\noindent No null Cartan curve $C$ in $R_{1}^{4}$ is a Bertrand curve if $%
\left\{ c^{\left( i\right) }\right\} _{1\leq i\leq 4}$ is linearly
independent.
\end{th}

\noindent {\it Proof} :\ {\bf \ }Let $C$\ be a Bertrand curve in $R_{1}^{4}$%
\ and $\overline{C}$ a Bertrand mate of $C$. $\overline{C}$ is distinct from
$C$. Let the pair $c\left( s\right) $ and $c\left( \overline{s}\right) $ be
of corresponding points of $C$ and $\overline{C}$ respectively. Then the
curve $\overline{C}$ is given by

\begin{equation}
\overline{c}\left( \overline{s}\right) =c\left( s\right) +\alpha \left(
s\right) W_{1}\left( s\right) ,  \label{a1}
\end{equation}%
where $\alpha $ is a function of $s$. Differentiating (\ref{a1}) with
respect to $s$, we obtain

\[
\frac{d\overline{s}}{ds}\frac{d\overline{c}\left( \overline{s}\right) }{d%
\overline{s}}=c^{\prime }\left( s\right) +\alpha ^{\prime }\left( s\right)
W_{1}\left( s\right) +\alpha \left( s\right) W_{1}^{\prime }\left( s\right)
.
\]%
Here the prime denotes the derivative with respect to $s$. By the Frenet
equations, it holds that

\[
\frac{d\overline{s}}{ds}\overline{L}\left( \overline{s}\right) =\left(
1-\alpha \left( s\right) k_{1}\left( s\right) \right) L\left( s\right)
-\alpha \left( s\right) N\left( s\right) +\alpha ^{\prime }\left( s\right)
W_{1}\left( s\right) .
\]%
Since $\left\langle \overline{L}\left( \overline{s}\right) ,\overline{W_{1}}%
\left( \overline{s}\right) \right\rangle =0$ and $\overline{W_{1}}\left(
\overline{s}\right) =\mp W_{1}\left( s\right) $, we obtain $\alpha ^{\prime
}\left( s\right) =0$, that is , $\alpha $ is a constant function with value $%
\alpha $ (we can use the same letter without confusion). Thus (\ref{a1}) is
rewritten as

\begin{equation}
\overline{c}\left( \overline{s}\right) =c\left( s\right) +\alpha W_{1}\left(
s\right) .  \label{a2}
\end{equation}%
Differentiating (\ref{a2}) with respect to $s$, we obtain

\begin{equation}
\frac{d\overline{s}}{ds}\overline{L}\left( \overline{s}\right) =\left(
1-\alpha k_{1}\left( s\right) \right) L\left( s\right) -\alpha N\left(
s\right) .  \label{a3}
\end{equation}%
From (\ref{a3}), we find

\[
0=-2\alpha \left( 1-\alpha k_{1}\left( s\right) \right) .
\]%
Then we get

\[
k_{1}\left( s\right) =\frac{1}{\alpha }{.}
\]%
So equation (\ref{a3}) becomes

\begin{equation}
\frac{d\overline{s}}{ds}\overline{L}\left( \overline{s}\right) =-\alpha
N\left( s\right) .  \label{a4}
\end{equation}%
Differentiating both sides of (\ref{a4}) with respect to $s$, we arrive at

\begin{equation}
\frac{d^{2}\overline{s}}{ds^{2}}\overline{L}\left( \overline{s}\right)
+\left( \frac{d\overline{s}}{ds}\right) ^{2}\overline{W_{1}}\left( \overline{%
s}\right) =-\alpha \left[ k_{1}\left( s\right) W_{1}\left( s\right)
+k_{2}\left( s\right) W_{2}\left( s\right) \right] .  \label{a5}
\end{equation}%
Since $\left\langle W_{2}\left( s\right) ,\overline{W_{1}}\left( \overline{s}%
\right) \right\rangle =0$, $\left\langle W_{2}\left( s\right) ,W_{1}\left(
s\right) \right\rangle =0$ and$\ \left\langle W_{2}\left( s\right) ,%
\overline{L}\left( \overline{s}\right) \right\rangle =0$, we get

\begin{equation}
-\alpha k_{2}\left( s\right) =0.\   \label{a6}
\end{equation}%
From the assumption in the hypothesis of the theorem, we have $k_{2}\left(
s\right) \neq 0$. So from (\ref{a6}), we obtain $\alpha =0$. Therefore (\ref%
{a2}) implies that $\overline{C}$ coincides with $C$. This is a
contradiction. This completes the proof.

\section{(1,2)-Bertrand curves in $R_{1}^{4}$}

Let $C$ be a null Cartan curve in $R_{1}^{4}$. We call $W_{1}$ the spacelike
Cartan 1-normal vector field along $C$, and the spacelike Cartan 1-normal
line of $C$ at $c\left( s\right) $ is a line generated by $W_{1}\left(
s\right) $ through $c\left( s\right) $. The spacelike Cartan (1,2)-normal
plane of $C$ at $c\left( s\right) $ is a plane spanned by $W_{1}\left(
s\right) $ and $W_{2}\left( s\right) $ through $c\left( s\right) $.

Let $C$ and $\overline{C}$ be null Cartan curves in $R_{1}^{4}$ and $\varphi
:s\rightarrow \overline{s}$ be a map such that each point $c\left( s\right) $
of $C$ corresponds to the point $\overline{c}\left( \overline{s}\right) =%
\overline{c}\left( \varphi \left( s\right) \right) $ of $\overline{C}.$ Here
$s$ and $\overline{s}$ are pseudo-arc parameters of $C$ and $\overline{C}$
respectively. If the Cartan (1,2)-normal plane at each point $c\left(
s\right) $ of $C$\ coincides with the Cartan (1,2)-normal plane at
corresponding point $\overline{c}\left( \overline{s}\right) =\overline{c}%
\left( \varphi \left( s\right) \right) $ of $\overline{C}$ and $\frac{%
d\varphi \left( s\right) }{ds}\neq 0$, then $C$ is called (1,2)-Bertrand
curve in $R_{1}^{4}$ and $\overline{C}$ is called (1,2)-Bertrand mate of $C$.

\begin{th}
Let $C$ be a null Cartan curve in $R_{1}^{4}$ with curvature functions $%
k_{1} $, $k_{2}$. Then $C$ is a (1,2)-Bertrand curve iff there are constants
$\alpha $ and $\beta $ satisfying either%
\[
{(i) \ \ \ \ \ \ \ }\alpha =0{ \ }and{ \ }\beta k_{2}\left( s\right)
\neq 1
\]%
or%
\[
{(ii) \ \ \ \ \ \ \ }\alpha \neq 0{ \ }and{ \ }\alpha k_{1}\left(
s\right) +\beta k_{2}\left( s\right) =1.
\]
\end{th}

\noindent {\it Proof}{\rm \ }: \ Now assume that $C$ is a null Cartan
(1,2)-Bertrand curve parametrized by the pseudo-arc parameter $s$ in $%
R_{1}^{4}$. The (1,2)-Bertrand mate $\overline{C}$ is given by

\begin{equation}
\overline{c}\left( \overline{s}\right) =c\left( s\right) +\alpha \left(
s\right) W_{1}\left( s\right) +\beta \left( s\right) W_{2}\left( s\right) .
\label{b1}
\end{equation}%
Here $\alpha $ and $\beta $ are functions and $\overline{s}$ is the
pseudo-arc parameter of $\overline{C}$. Since the plane spanned by $W_{1}$
and $W_{2}$ coincides with the plane spanned by $\overline{W_{1}}$ and $%
\overline{W_{2}}$, we can write

\begin{eqnarray}
\overline{W_{1}}\left( \overline{s}\right) &=&\cos \theta \left( s\right)
W_{1}\left( s\right) +\sin \theta \left( s\right) W_{2}\left( s\right)
\label{b2} \\
\overline{W_{2}}\left( \overline{s}\right) &=&-\sin \theta \left( s\right)
W_{1}\left( s\right) +\cos \theta \left( s\right) W_{2}\left( s\right)
\nonumber
\end{eqnarray}%
If we differentiate (\ref{b1}) with respect to $s$, then we get

\begin{eqnarray}
\overline{L}\left( \overline{s}\right) \frac{d\overline{s}}{ds} &=&\left[
1-\alpha \left( s\right) k_{1}\left( s\right) -\beta \left( s\right)
k_{2}\left( s\right) \right] L\left( s\right) -\alpha \left( s\right)
N\left( s\right)  \label{b3} \\
&&+\alpha ^{\prime }\left( s\right) W_{1}\left( s\right) +\beta ^{\prime
}\left( s\right) W_{2}\left( s\right) .  \nonumber
\end{eqnarray}%
If we use (\ref{b2}) and (\ref{b3})\ together, then we have

\begin{eqnarray*}
0 &=&\left\langle \overline{W_{1}}\left( \overline{s}\right) ,\overline{L}%
\left( \overline{s}\right) \frac{d\overline{s}}{ds}\right\rangle =\alpha
^{\prime }\left( s\right) \cos \theta \left( s\right) +\beta ^{\prime
}\left( s\right) \sin \theta \left( s\right) \\
0 &=&\left\langle \overline{W_{2}}\left( \overline{s}\right) ,\overline{L}%
\left( \overline{s}\right) \frac{d\overline{s}}{ds}\right\rangle =-\alpha
^{\prime }\left( s\right) \sin \theta \left( s\right) +\beta ^{\prime
}\left( s\right) \cos \theta \left( s\right) .
\end{eqnarray*}%
So we obtain $\ \alpha ^{\prime }\left( s\right) =\beta ^{\prime }\left(
s\right) =0$, that is $\alpha $ and $\beta $ are constants. Then we can
rewrite (\ref{b3}) as

\begin{equation}
\overline{L}\left( \overline{s}\right) \frac{d\overline{s}}{ds}=\left[
1-\alpha k_{1}\left( s\right) -\beta k_{2}\left( s\right) \right] L\left(
s\right) -\alpha N\left( s\right) .  \label{b4}
\end{equation}%
From (\ref{b4}), we get

\[
0=\left\langle \overline{L}\left( \overline{s}\right) \frac{d\overline{s}}{ds%
},\overline{L}\left( \overline{s}\right) \frac{d\overline{s}}{ds}%
\right\rangle =-2\alpha \left[ 1-\alpha k_{1}\left( s\right) -\beta
k_{2}\left( s\right) \right] {.}
\]%
Here observe that it must be either

\[
{(i) \ \ \ \ \ \ \ \ \ }\alpha =0\,\,\,\,{ and }\,\,\,\,\beta
k_{2}\left( s\right) \neq 1
\]%
or

\[
{(ii) \ \ \ \ \ \ \ \ \ }\alpha \neq 0\,\,\,\,{ and }\,\,\,\,\alpha
k_{1}\left( s\right) +\beta k_{2}\left( s\right) )=1.
\]%
Note that\ if\ both\ $\alpha $ and $\left[ 1-\alpha \left( s\right)
k_{1}\left( s\right) -\beta \left( s\right) k_{2}\left( s\right) \right] $
are zero in equation (\ref{b4}), this leads us to conclude that $\frac{d%
\overline{s}}{ds}=0$. But we know from the definition of (1,2)-Bertrand
curve the map $\varphi \left( s\right) =\overline{s}$ under which the points
$c\left( s\right) $ of $C$ and $\overline{c}\left( \overline{s}\right) $ of $%
\overline{C}$ correspond is a regular map. Therefore the case $\frac{d%
\overline{s}}{ds}=0$ in (\ref{b4}) never happens. To characterize a
(1,2)-Bertrand curve in $R_{1}^{4}$, we delete the case $\frac{d\overline{s}%
}{ds}=\frac{d\varphi \left( s\right) }{ds}=0$, and always think the position
$\frac{d\overline{s}}{ds}\neq 0$.

Now we assume the contrary cases and prove that the curve $C$ is a
(1,2)-Bertrand curve.

(i) {\bf \ }$\alpha =0${\bf \ \ }and{\bf \ \ }$\beta k_{2}\left( s\right)
)\neq 1${\bf .}

Assume that $C$ is a null Cartan curve in $R_{1}^{4}$ whose curvature
functions$\ k_{1}\left( s\right) $ and $k_{2}\left( s\right) $ satisfy the
relation\ $\beta k_{2}\left( s\right) )\neq 1$ for some constant real
numbers $\alpha =0$ and $\beta $. Now we define a curve $\overline{C}$ by

\begin{equation}
\overline{c}\left( s\right) =c\left( s\right) +\alpha W_{1}\left( s\right)
+\beta W_{2}\left( s\right) .  \label{b5}
\end{equation}%
where $s$ is the pseudo-arc parameter of $C$. Differentiating (\ref{b5})
with respect to $s$ and using the Frenet equations, we obtain

\[
\frac{d\overline{c}\left( s\right) }{ds}=\left[ 1-\alpha k_{1}\left(
s\right) -\beta k_{2}\left( s\right) \right] L\left( s\right) -\alpha
N\left( s\right) .
\]%
By using (i),we get

\begin{equation}
\frac{d\overline{c}\left( s\right) }{ds}=\left[ 1-\beta k_{2}\left( s\right) %
\right] L\left( s\right) .  \label{b6}
\end{equation}%
Because of (i), the curve $\overline{C}$ is a regular curve. Let us define a
regular map $\varphi :s\rightarrow \overline{s}$ by

\[
\overline{s}=\varphi \left( s\right) =\int\limits_{0}^{s}\left\langle \frac{%
d^{2}\overline{c}\left( s\right) }{ds^{2}},\frac{d^{2}\overline{c}\left(
s\right) }{ds^{2}}\right\rangle ^{\frac{1}{4}}ds,
\]%
where $\overline{s}$ denotes the pseudo-arc parameter of $\overline{C}$.
Then we obtain

\begin{equation}
\frac{d\overline{s}}{ds}=\frac{d\varphi \left( s\right) }{ds}=\sqrt{%
\left\vert 1-\beta k_{2}\left( s\right) \right\vert }>0.  \label{b7}
\end{equation}%
Thus the curve $\overline{C}$ is rewritten as

\begin{equation}
\overline{c}\left( \overline{s}\right) =\overline{c}\left( \varphi \left(
s\right) \right) =c\left( s\right) +\beta W_{2}\left( s\right) .  \label{b8}
\end{equation}%
If we differentiate (\ref{b8}) and use the Frenet equations for the null
Cartan curve in $R_{1}^{4}$, we have

\begin{equation}
\overline{L}\left( \overline{s}\right) \frac{d\overline{s}}{ds}=(1-\beta
k_{2}\left( s\right) )L\left( s\right) .  \label{b9}
\end{equation}%
We find the vector $\overline{N}$ uniquely by using the formula in \cite%
{duggal}. We get

\begin{equation}
\overline{N}\left( \overline{s}\right) (1-\beta k_{2}\left( s\right) )=\frac{%
d\overline{s}}{ds}N\left( s\right) .  \label{b10}
\end{equation}%
If we differentiate (\ref{b9}), then we obtain

\begin{equation}
\overline{W}_{1}\left( \overline{s}\right) \left( \frac{d\overline{s}}{ds}%
\right) ^{2}+\overline{L}\left( \overline{s}\right) \frac{d^{2}\left(
\overline{s}\right) }{ds^{2}}=(1-\beta k_{2}\left( s\right) )^{\prime
}L\left( s\right) +(1-\beta k_{2}\left( s\right) )W_{1}\left( s\right) .
\label{b11}
\end{equation}%
From \ (\ref{b11}) we have

\[
\left( \frac{d\overline{s}}{ds}\right) ^{2}=\left\vert 1-\beta k_{2}\left(
s\right) \right\vert \neq 0.
\]%
Let us take $\left( \frac{d\overline{s}}{ds}\right) ^{2}=1-\beta k_{2}\left(
s\right) $. Then equations (\ref{b9}) and\ (\ref{b10}) become

\begin{equation}
\overline{L}\left( \overline{s}\right) =\frac{d\overline{s}}{ds}L\left(
s\right) ,  \label{b12}
\end{equation}

\begin{equation}
\overline{N}\left( \overline{s}\right) \frac{d\overline{s}}{ds}=N\left(
s\right) .  \label{b13}
\end{equation}%
Now differentiating$\ $(\ref{b12}) with respect to $s$ and using Frenet
equations, we get

\begin{equation}
\overline{W}_{1}\left( \overline{s}\right) \frac{d\overline{s}}{ds}=\frac{%
d^{2}\left( \overline{s}\right) }{ds^{2}}L\left( s\right) +\frac{d\overline{s%
}}{ds}W_{1}\left( s\right) .  \label{b14}
\end{equation}%
Using the equations (\ref{b13}) and (\ref{b14}) together, we obtain

\[
0=\left\langle \overline{N}\left( \overline{s}\right) \frac{d\overline{s}}{ds%
},\overline{W}_{1}\left( \overline{s}\right) \frac{d\overline{s}}{ds}%
\right\rangle =\frac{d^{2}\left( \overline{s}\right) }{ds^{2}}.
\]%
From here we can say that $\frac{d\overline{s}}{ds}=\ell _{0}$ is a
constant. Then the equations (\ref{b12}) and (\ref{b13}) shrink to

\begin{eqnarray}
\overline{L}\left( \overline{s}\right) &=&\ell _{0}L\left( s\right) ,
\label{b15} \\
\overline{N}\left( \overline{s}\right) &=&\frac{1}{\ell _{0}}N\left(
s\right) .  \label{b16}
\end{eqnarray}%
If we differentiate (\ref{b15}), then we get

\[
\overline{W}_{1}\left( \overline{s}\right) \frac{d\overline{s}}{ds}=\ell
_{0}W_{1}\left( s\right) .
\]%
Using the fact that $\frac{d\overline{s}}{ds}=\ell _{0}$, we have

\begin{equation}
\overline{W}_{1}\left( \overline{s}\right) =W_{1}\left( s\right) .
\label{b17}
\end{equation}%
If we differentiate (\ref{b16}), we have

\begin{equation}
\left[ \overline{k}_{1}\left( \overline{s}\right) \overline{W}_{1}\left(
\overline{s}\right) +\overline{k}_{2}\left( \overline{s}\right) \overline{W}%
_{2}\left( \overline{s}\right) \right] \ell _{0}=\frac{1}{\ell _{0}}\left[
k_{1}\left( s\right) W_{1}\left( s\right) +k_{2}\left( s\right) W_{2}\left(
s\right) \right] .  \label{x1}
\end{equation}%
From (\ref{x1}),we get

\begin{equation}
\left[ \left( \overline{k}_{1}\left( \overline{s}\right) \right) ^{2}+\left(
\overline{k}_{2}\left( \overline{s}\right) \right) ^{2}\right] \left( \ell
_{0}\right) ^{4}=\left[ \left( k_{1}\left( s\right) \right) ^{2}+\left(
k_{2}\left( s\right) \right) ^{2}\right] .  \label{b18}
\end{equation}%
Differentiating (\ref{b17}), we obtain

\begin{equation}
\left[ -\overline{k}_{1}\left( \overline{s}\right) \overline{L}\left(
\overline{s}\right) -\overline{N}\left( \overline{s}\right) \right] \ell
_{0}=-k_{1}\left( s\right) L\left( s\right) -N\left( s\right) .  \label{x2}
\end{equation}%
By using (\ref{x2}), we get

\begin{equation}
\overline{k}_{1}\left( \overline{s}\right) =\frac{k_{1}\left( s\right) }{%
\left( \ell _{0}\right) ^{2}}.  \label{b19}
\end{equation}%
Now putting (\ref{b19}) into the equation (\ref{b18}), we have

\[
\left\vert \overline{k}_{2}\left( \overline{s}\right) \right\vert
=\left\vert \frac{k_{2}\left( s\right) }{\left( \ell _{0}\right) ^{2}}%
\right\vert .
\]
Let us take

\begin{equation}
\overline{k}_{2}\left( \overline{s}\right) =\frac{k_{2}\left( s\right) }{%
\left( \ell _{0}\right) ^{2}}.  \label{b20}
\end{equation}%
If we use (\ref{b17}), (\ref{b19}) and (\ref{b20}) into (\ref{x1}), we obtain

\begin{equation}
\overline{W}_{2}=W_{2}.  \label{b21}
\end{equation}%
Observe that the coincidence of the plane spanned by $W_{1}$ and $W_{2}$
with the plane spanned by $\overline{W}_{1}$ and \ $\overline{W}_{2}$ is
trivial. Therefore $C$ is a (1,2)-Bertrand curve in $R_{1}^{4}$.

\noindent (ii) $\alpha \neq 0$ \ and \ $\alpha k_{1}+\beta k_{2}=1$.

Now we assume that $C$ is a null Cartan curve whose curvature functions
satisfy the equation $\ \alpha k_{1}+\beta k_{2}=1$\ for some constant real
numbers $\alpha \neq 0$ and $\beta $. Think about a curve $\overline{C}$
defined as

\begin{equation}
\overline{c}\left( s\right) =c\left( s\right) +\alpha W_{1}\left( s\right)
+\beta W_{2}\left( s\right) .  \label{b22}
\end{equation}%
Differentiating (\ref{b22}) with respect to $s$ and using the Frenet
equations, we get

\[
\frac{d\overline{c}\left( s\right) }{ds}=\left[ 1-\alpha k_{1}\left(
s\right) -\beta k_{2}\left( s\right) \right] L\left( s\right) -\alpha
N\left( s\right) .
\]%
By using (ii), we have

\begin{equation}
\frac{d\overline{c}\left( s\right) }{ds}=-\alpha N\left( s\right) .
\label{b23}
\end{equation}%
Because of (ii), the curve $\overline{C}$ is a regular curve. Then we define
a regular map $\varphi :s\rightarrow \overline{s}$ by

\[
\overline{s}=\varphi \left( s\right) =\int\limits_{0}^{s}\left\langle \frac{%
d^{2}\overline{c}\left( s\right) }{ds^{2}},\frac{d^{2}\overline{c}\left(
s\right) }{ds^{2}}\right\rangle ^{\frac{1}{4}}ds,
\]%
where $\overline{s}$ denotes the pseudo-arc parameter of $\overline{C}.$
Then we obtain

\begin{equation}
\frac{d\overline{s}}{ds}=\frac{d\varphi \left( s\right) }{ds}=\left[ \alpha
^{2}\left( \left( k_{1}\right) ^{2}+\left( k_{2}\right) ^{2}\right) \right]
^{\frac{1}{4}}>0.  \label{b24}
\end{equation}%
Now we can rewrite equation (\ref{b22}) as%
\[
\overline{c}\left( \overline{s}\right) =c\left( s\right) +\alpha W_{1}\left(
s\right) +\beta W_{2}\left( s\right) .
\]%
Differentiating this equation with respect to $s$, we get

\[
\frac{d\overline{s}}{ds}\overline{L}\left( \overline{s}\right) =\left[
1-\alpha k_{1}\left( s\right) -\beta k_{2}\left( s\right) \right] L\left(
s\right) -\alpha N\left( s\right) .
\]%
Inserting (ii) in the above, we have

\begin{equation}
\frac{d\overline{s}}{ds}\overline{L}\left( \overline{s}\right) =-\alpha
N\left( s\right) .  \label{b26}
\end{equation}%
Once again we use the formula in \cite{duggal} to find $\overline{N}.$ We
find $\overline{N}$\ uniquely as

\begin{equation}
\overline{N}\left( \overline{s}\right) =-\frac{d\overline{s}}{ds}\frac{1}{%
\alpha }L\left( s\right)  \label{b27}
\end{equation}%
Differentiating (\ref{b26}), we obtain

\begin{equation}
\frac{d^{2}\overline{s}}{ds^{2}}\overline{L}\left( \overline{s}\right)
+\left( \frac{d\overline{s}}{ds}\right) ^{2}\overline{W}_{1}\left( \overline{%
s}\right) =-\alpha \left[ k_{1}\left( s\right) W_{1}\left( s\right)
+k_{2}\left( s\right) W_{2}\left( s\right) \right] .  \label{b28}
\end{equation}%
Using (\ref{b27}) and (\ref{b28}), we get

\[
\left\langle \overline{N}\left( \overline{s}\right) ,\left[ \frac{d^{2}%
\overline{s}}{ds^{2}}\overline{L}\left( \overline{s}\right) +\left( \frac{d%
\overline{s}}{ds}\right) ^{2}\overline{W}_{1}\left( \overline{s}\right) %
\right] \right\rangle
\]

\[
=\left\langle -\frac{d\overline{s}}{ds}\frac{1}{\alpha }L\left( s\right) ,%
\left[ -\alpha k_{1}\left( s\right) W_{1}\left( s\right) -\alpha k_{2}\left(
s\right) W_{2}\left( s\right) \right] \right\rangle =0.
\]%
From the above fact, we get

\[
\frac{d^{2}\overline{s}}{ds^{2}}=0.
\]%
So $\frac{d\overline{s}}{ds}=\ell _{0}$ is a constant. Using this in (\ref%
{b28}), we have

\begin{equation}
\overline{W}_{1}\left( \overline{s}\right) =-\frac{\alpha k_{1}}{\left( \ell
_{0}\right) ^{2}}W_{1}\left( s\right) -\frac{\alpha k_{2}}{\left( \ell
_{0}\right) ^{2}}W_{2}\left( s\right) .  \label{b29}
\end{equation}%
By using (\ref{b29}), we obtain

\begin{equation}
\left( \ell _{0}\right) ^{4}=\alpha ^{2}\left[ \left( k_{1}\right)
^{2}+\left( k_{2}\right) ^{2}\right] .  \label{b30}
\end{equation}%
So we can write

\begin{equation}
\overline{W}_{1}\left( \overline{s}\right) =\cos \tau \left( s\right)
W_{1}\left( s\right) +\sin \tau \left( s\right) W_{2}\left( s\right) .
\label{b31}
\end{equation}%
Differentiating (\ref{b31}), we arrive at

\begin{eqnarray*}
\left[ -\overline{k}_{1}\left( \overline{s}\right) \overline{L}\left(
\overline{s}\right) -\overline{N}\left( \overline{s}\right) \right] \ell
_{0} &=&\frac{d}{ds}\left[ \cos \tau \left( s\right) \right] W_{1}\left(
s\right) \\
&&+\cos \tau \left( s\right) \left[ -k_{1}\left( s\right) L\left( s\right)
-N\left( s\right) \right] \\
&&+\frac{d}{ds}\left[ \sin \tau \left( s\right) \right] W_{2}\left( s\right)
\\
&&+\sin \tau \left( s\right) \left[ -k_{2}\left( s\right) L\left( s\right) %
\right] .
\end{eqnarray*}%
Using (\ref{b26}) and (\ref{b27}) in the above, we have

\begin{eqnarray}
\left[ -\overline{k}_{1}\left( \overline{s}\right) \left( -\frac{\alpha }{%
\ell _{0}}N\left( s\right) \right) -\left( -\frac{\ell _{0}}{\alpha }L\left(
s\right) \right) \right] \ell _{0} &=&\frac{d}{ds}\left[ \cos \tau \left(
s\right) \right] W_{1}\left( s\right)  \label{b32} \\
&&+\cos \tau \left( s\right) \left[ -k_{1}\left( s\right) L\left( s\right)
-N\left( s\right) \right]  \nonumber \\
&&+\frac{d}{ds}\left[ \sin \tau \left( s\right) \right] W_{2}\left( s\right)
\nonumber \\
&&+\sin \tau \left( s\right) \left[ -k_{2}\left( s\right) L\left( s\right) %
\right] .  \nonumber
\end{eqnarray}%
Using the fact $\left\langle W_{1}\left( s\right) ,L\left( s\right)
\right\rangle =\left\langle W_{1}\left( s\right) ,N\left( s\right)
\right\rangle =\left\langle W_{2}\left( s\right) ,L\left( s\right)
\right\rangle =\left\langle W_{2}\left( s\right) ,N\left( s\right)
\right\rangle =0$ in (\ref{b32}), we obtain

\[
\frac{d}{ds}\left[ \cos \tau \left( s\right) \right] =\frac{d}{ds}\left[
\sin \tau \left( s\right) \right] =0.
\]%
So the function $\tau \left( s\right) $ must be the constant $\tau \left(
s_{0}\right) $. Using (\ref{b32}), we can also get

\[
\overline{k}_{1}\left( \overline{s}\right) \alpha =-\cos \tau \left(
s\right) =\frac{\alpha k_{1}\left( s\right) }{\left( \ell _{0}\right) ^{2}}%
{.}
\]%
And from the above, we have

\begin{equation}
\overline{k}_{1}\left( \overline{s}\right) =\frac{k_{1}\left( s\right) }{%
\left( \ell _{0}\right) ^{2}}{.}  \label{b33}
\end{equation}%
Differentiating (\ref{b27}), we obtain

\begin{equation}
\left[ \overline{k}_{1}\left( \overline{s}\right) \overline{W}_{1}\left(
\overline{s}\right) +\overline{k}_{2}\left( \overline{s}\right) \overline{W}%
_{2}\left( \overline{s}\right) \right] \ell _{0}=-\frac{\ell _{0}}{\alpha }%
W_{1}\left( s\right) .  \label{b34}
\end{equation}%
From (\ref{b34}), we get

\begin{equation}
(\overline{k}_{1}\left( \overline{s}\right) )^{2}+\left( \overline{k}%
_{2}\left( \overline{s}\right) \right) ^{2}=\frac{1}{\alpha ^{2}}.
\label{b35}
\end{equation}%
Using (\ref{b30}) and (\ref{b33}) into the equation (\ref{b35}), we get

\[
\left\vert \overline{k}_{2}\left( \overline{s}\right) \right\vert
=\left\vert \frac{k_{2}\left( s\right) }{(\ell _{0})^{2}}\right\vert .
\]%
Choosing,

\[
\overline{k}_{2}\left( \overline{s}\right) =-\frac{k_{2}\left( s\right) }{%
(\ell _{0})^{2}},
\]%
and putting it in (\ref{b34}), we obtain

\[
\left\{ \frac{k_{1}\left( s\right) }{\left( \ell _{0}\right) ^{2}}\left( -%
\frac{\alpha k_{1}\left( s\right) }{\left( \ell _{0}\right) ^{2}}W_{1}\left(
s\right) -\frac{\alpha k_{2}\left( s\right) }{(\ell _{0})^{2}}W_{2}\left(
s\right) \right) +\left( -\frac{k_{2}\left( s\right) }{(\ell _{0})^{2}}%
\right) \overline{W}_{2}\left( \overline{s}\right) \right\} \ell _{0}=-\frac{%
\ell _{0}}{\alpha }W_{1}\left( s\right) .
\]%
Simplifying this equation, we arrive at

\[
\overline{W}_{2}\left( \overline{s}\right) =\frac{\alpha k_{2}\left(
s\right) }{(\ell _{0})^{2}}W_{1}\left( s\right) -\frac{\alpha k_{1}\left(
s\right) }{\left( \ell _{0}\right) ^{2}}W_{2}\left( s\right) =-\sin \tau
\left( s_{0}\right) W_{1}\left( s\right) +\cos \tau \left( s_{0}\right)
W_{2}\left( s\right) .
\]%
And it is trivial that the Cartan (1,2)-normal plane at each point $c\left(
s\right) $ of $C$ coincides with the Cartan (1,2)-normal plane at
corresponding point $\overline{c}\left( \overline{s}\right) $ of $\overline{C%
}$. Therefore $C$ is a (1,2)-Bertrand curve in $R_{1}^{4}$.

\section{An example of (1,2)-Bertrand curve}

Let $a$, $b$ be nonzero contants such that $a\neq \mp b$, and let $\ C$ be
the curve in\ $R_{1}^{4}$ defined by%
\[
c\left( s\right) =\frac{1}{\sqrt{a^{2}+b^{2}}}\left[ \frac{1}{a}\sinh as,%
\frac{1}{a}\cosh as,\frac{1}{b}\sin bs,\frac{1}{b}\cos bs\right] .
\]%
Then we get the Cartan frame and the curvature\ functions as follows:

\begin{eqnarray*}
L\left( s\right) &=&\frac{1}{\sqrt{a^{2}+b^{2}}}\left[ \cosh as,\sinh
as,\cos bs,-\sin bs\right] , \\
W_{1}\left( s\right) &=&\frac{1}{\sqrt{a^{2}+b^{2}}}\left[ a\sinh as,a\cosh
as,-b\sin bs,-b\cos bs\right] , \\
N\left( s\right) &=&-\frac{\sqrt{a^{2}+b^{2}}}{2}\left[ \cosh as,\sinh
as,-\cos bs,\sin bs\right] , \\
W_{2}\left( s\right) &=&\frac{1}{\sqrt{a^{2}+b^{2}}}\left[ b\sinh as,b\cosh
as,a\sin bs,a\cos bs\right] , \\
k_{1} &=&\frac{b^{2}-a^{2}}{2}, \\
k_{2} &=&-ab.
\end{eqnarray*}%
To enlighten the case (i), we take the constants real numbers $\alpha $ and $%
\beta $ as follows:

\[
\alpha =0{ \ and \ }\beta =\frac{1}{ab}.
\]%
Then

\[
\alpha k_{1}+\beta k_{2}=-1\neq 1
\]%
holds. Therefore $C$ is a (1,2)-Bertrand curve in $R_{1}^{4}$, and its
Bertrand mate $\overline{C}$ is given by

\[
\overline{c}\left( \overline{s}\right) =\frac{2}{\sqrt{a^{2}+b^{2}}}\left[
\frac{1}{a}\sinh \left( a\frac{\overline{s}}{\sqrt{2}}\right) ,\frac{1}{a}%
\cosh \left( a\frac{\overline{s}}{\sqrt{2}}\right) ,\frac{1}{b}\sin \left( b%
\frac{\overline{s}}{\sqrt{2}}\right) ,\frac{1}{b}\cos \left( b\frac{%
\overline{s}}{\sqrt{2}}\right) \right] ,
\]%
where $\overline{s}$ is the pseudo-arc parameter of $\overline{C}$, and a
regular map $\varphi :s\rightarrow \overline{s}$ is given by

\[
\overline{s}=\varphi \left( s\right) =\sqrt{2}s.
\]%
For the case (ii), we take real constants $\alpha $ and $\beta $ as follows:

\[
\alpha =\frac{1}{b^{2}-a^{2}}{ and }\beta =-\frac{1}{2ab}.
\]%
Then it is trivial the following equation

\[
\alpha k_{1}+\beta k_{2}=1
\]%
holds. Therefore $C$ is a (1,2)-Bertrand curve in $R_{1}^{4}$, and its
Bertrand mate is given by

\[
\overline{c}\left( \overline{s}\right) =\frac{(\ell _{0})^{2}}{\sqrt{%
a^{2}+b^{2}}}\left[ \frac{1}{a}\sinh \left( \frac{a}{\ell _{0}}\overline{s}%
\right) ,\frac{1}{a}\cosh \left( \frac{a}{\ell _{0}}\overline{s}\right) ,-%
\frac{1}{b}\sin \left( \frac{b}{\ell _{0}}\overline{s}\right) ,-\frac{1}{b}%
\cos \left( \frac{b}{\ell _{0}}\overline{s}\right) \right] .
\]%
Here $\overline{s}$ is the pseudo-arc parameter of $\overline{C}$, $\frac{d%
\overline{s}}{ds}=\ell _{0}$, and a regular map $\varphi :s\rightarrow
\overline{s}$ is given by

\[
\overline{s}=\varphi \left( s\right) =\sqrt{\frac{a^{2}+b^{2}}{2\left(
b^{2}-a^{2}\right) }}s.
\]

{\ \noindent Mehmet G\"{o}\c{c}men}

{\ \noindent Department of Mathematics, Faculty of Arts and Sciences, \.{I}n%
\"{o}n\"{u} University }

{\ \noindent 44280 Malatya, Turkey }

{\ \noindent Email: mgocmen1903@gmail.com}

{\ \noindent Sad\i k Kele\c{s} }

{\ \noindent Department of Mathematics, Faculty of Arts and Sciences, \.{I}n%
\"{o}n\"{u} University }

{\ \noindent 44280 Malatya, Turkey }

{\ \noindent Email: skeles@inonu.edu.tr }

\end{document}